\documentclass[11pt,epsf]{amsart}
\usepackage{epsfig}
\usepackage{amssymb}
\usepackage{eepic}
\usepackage{mathrsfs}
\usepackage{xypic}
\usepackage{verbatim}
\usepackage{graphicx}
\usepackage{psfrag}
\CompileMatrices

\def\makemargins{
	\oddsidemargin .25in
	\evensidemargin .25in
	\textwidth 6.00in
}
\makemargins
\theoremstyle{plain}
\newtheorem{theorem}[subsection]{Theorem}
\newtheorem{lemma}[subsection]{Lemma}
\newtheorem{proposition}[subsection]{Proposition}
\newtheorem{corollary}[subsection]{Corollary}

\theoremstyle{definition}
\newtheorem{definition}[subsection]{Definition}
\newtheorem{example}[subsection]{Example}

\theoremstyle{remark}

\newcount\timehh\newcount\timemm
\timehh=\time 
\divide\timehh by 60 \timemm=\time
\newcommand{\draftauthor}[1]{\author{#1
    {
      --- \protect \protect\sc\today\ ---
      \ifnum\timehh<10 0\fi\number\timehh\,:\,\ifnum\timemm<10 0\fi\number\timemm
      \protect \, \, \protect \bf DRAFT
    }
  }
}
\begin{document}
\title{Reconstructing Metric Trees from Order Information 
on Triples is NP Complete}
\author{ Eric K. Babson}
%\address{Department of Mathematics\\
%University of Washington\\
%Seattle, WA}
\email{babson@math.washington.edu}

\maketitle

\begin{abstract}
We show that reconstructing a tree from order information 
on triples is NP-hard.  This is in contrast to the case 
for ultra-metrics and for subtree information on quadruples 
which are both known to allow polynomial time reconstruction.  
\end{abstract}

\begin{section}{Introduction}
This paper deals with the computational complexity of finding a tree 
compatible with known combinatorial data.  
This type of problem arises among other places \cite{b} 
in trying to reconstruct a 
phylogenetic tree from partial information.  
We might for example know about all triples of species 
which two are most closely related and which two least.  
Equivalently we might know for each pair of species which of the 
others are more closely related to the first and which to the second.  
To get decision problems these types of data are abstracted to 
triples structures and midpoints structures respectively.  
We then address two problems.  We show that it is computationally
difficult (NP-complete and NP-hard respectively) to either determine 
whether the known data is compatible with any tree structure or given 
that it is compatible with some tree structure to find such a tree.  
A computationally equivalent formulation  of a tree structure is 
the information about each quadruple of species of the tree structure 
underlying just these four \cite{g7}.

Triples and midpoints structures are computationally closely related:  
There is a polynomial time bijection between them which preserves the set 
of compatible tree geometries.  
Thus we will focus entirely on midpoints structures but get the 
same results for both.  
There is a polynomial time algorithm to test whether a given tree 
geometry is compatible with a given midpoints 
structure so that determining whether there exists a tree 
compatible with a given structure is in NP.  

Some related questions have been previously studied.  
It is known \cite{g7} that if the metric tree is 
required to have some point in the tree equidistant from every element
(an ultra-metric) the computational problems can both be solved in polynomial 
time.  
A less natural restriction than an ultra-metric is that every edge of 
the tree into which the set is embedded contain the 
midpoint between some pair of elements.  
Call this tree geometry the midpoints tree.  
This idea arises from 
the fact that any edge containing the midpoint between two of the 
elements can be quickly identified from the above types of data 
\cite{bgp} and hence so can the midpoints tree.  
There is a metric on the midpoints tree studied in \cite{bgp}, 
\cite{g7} and \cite{g} where it is referred to as the order 
distance.  Unfortunately even in cases where 
the given midpoints structure has a realization with the midpoints tree 
geometry the order distance need not be such a realization since it need not 
be compatible with the original midpoints structure.   
In the ultra-metric situation the midpoints tree is the unique minimal 
tree geometry realizing any given midpoints structure which has an 
ultra-metric realization.  
Warnow has asked whether the midpoints tree is the unique minimal 
tree geometry realizing any given realizable midpoints structure in general;
see also Theorem 2 of \cite{g7}.  
The examples in this paper porvide counterexamples for the above.  

In section 2 we introduce notation.  Section 3 states the results.  
In section 4 we give the main construction: a standard NP-complete 
problem (3-SAT) is encoded by a midpoints structure.
In sections 5 and 6 we show the equivalence of the satisfiability of 
any case of 3-SAT with the realizability of the midpoints structure 
which encodes that case.  

\end{section}
\begin{section}{Notation}
For notational convenience only generic structures (for which none of the 
distances are equal) will be considered.  The results still hold for 
arbitrary data since every realizable triples or midpoints structure without 
ties can be realized by a tree without ties.  

It will be convenient throughout to fix a total order on the 
elements of the finite sets.  

\begin{definition}
A {\bf triples structure} on an ordered finite set $X$ 
is a set $\{<_x\}_{x \in X}$ of 
relations on $X$ such that for all $x,y,z \in X$ with $x \not= y$ 
we have $x\not<_z x$, $x<_xy$ and either $y <_zx$ or $x <_zy$ 
but not both.  
\end{definition}

A metric on $X$ will be a realization of a triples structure $\{<_x\}_{x \in X}$ if 
the distance from $x$ to $z$ is less than the distance from $y$ to $z$ 
whenever $x<_zy$.  

Write ${X\choose 2}$ for $\{A \subseteq X||A|=2\}$.  

\begin{definition}
A {\bf midpoints structure} on an ordered set $X$ is a 
map $m:{X\choose 2} \rightarrow 2^X$ with 
\{$\hbox{max}\{x,y\}\}= \{x,y\}\cap m\{x,y\}$ for all $x \not= y \in X$.  
\end{definition}

A metric on $X$ will be a realization of a midpoints structure $m$ if 
the distance from $x$ to $z$ is less than the distance from $y$ to $z$ 
whenever $y<x$ and $z \in m\{x,y\}$ or $x<y$ and $z \not\in m\{x,y\}$.  

There are polynomial time algorithms to translate between these two 
types of data: 

To get a midpoints structure $m$ from a triples structure $\{<_x\}_{x \in X}$: 
For every $x<y$ set $m\{x,y\}=\{ z \in X|y<_z x\}$.  

To get a triples structure $\{<_x\}_{x \in X}$ from a midpoints structure $m$: 
For every $y<x\in X$ and $z \in X$ if $z \in m\{x,y\}$ set $x<_zy$ and otherwise 
set $y<_zx$.  
\begin{lemma}
These algorithms are inverses to each other and a metric 
is a realization of  a given midpoints structure if and only if it is a realization of  the associated triples structure.  
\end{lemma}

If $S \subseteq X$ write $S^1=S$ and $S^c=X \backslash S$.  
Think of $S \subseteq X$ as an edge of a tree with leaves labeled by $X$ 
for which the edge splits $X$ into the sets $S^1$ and $S^c$.  

\begin{definition}\label{int}
If $S,T \subseteq X$ define $[S,T]=\{U\subseteq X|\hbox{ for some }e,f,g \in \{c, 1\}, S^e \subseteq U^g 
\subseteq T^f\}$, and define $[S,T)=[S,T] \backslash \{T\}$.  If $U \in [S,T]$ 
then write $S:U:T$.  
\end{definition}

Thus $[S,T]$ is all edges which appear on the path between the edges $S$ and $T$ in 
some tree.  
If $t:=\sum_{e \in E}t(e)x_e \in {\bf R}^E$ and $W \subseteq E$ then write 
$t_W:=\|\pi_W t \|_1=\Sigma_{e \in W}t(e)$.  
A metric tree will be represented as an element $t \in {\bf R}^{2^X}$ with 
$t_{\{S\}}=t(S)=t(S^c)$ being 
the length of the edge splitting the leaves into $S^1$ and $S^c$.  
Write $R_1:R_2: \ldots :R_r$ if $R_i:R_j:R_k$ for every $1\leq i<j<k\leq r$.  
Note that if any two of $R_1:R_2:R_3$, $R_1:R_2:R_4$, 
$R_1:R_3:R_4$ and $R_2:R_3:R_4$ hold then all four hold.  

A structure will be called realizable if there is a tree metric which is a realization of  
it.  More precisely, 

\begin{definition}\label{real}
A midpoints structure $m$ on $X$ is called 
{\bf realizable with realization $t \in {\bf R}^{2^X}$} 
and {\bf tree structure} $\{S \subseteq X|t(S)>0\}$ if 
\begin{enumerate}
\item if $S \subseteq X$ then $t(S)=t(S^c)\geq 0$,
\item if $S, T \subseteq X$ with $S^e \cap T^f \not= \emptyset$ 
for every $e,f \in \{1,c\}$ then $t(S)t(T)=0$ and
\item if $x\not= x' \in X$ then $t_{[\{x\},m\{x,x'\}]}>t_{[\{x'\},m\{x,x'\})}$.
\end{enumerate}
\end{definition}

The first condition ensures that the length $t(S)$ of the edge $S$ is the 
same as that of $S^c$ (which is the same edge) and is nonnegative.  
The second condition ensures that the nonzero edges form a tree.  
The third condition ensures that the edge $m\{x,x'\}$ contains in its interior 
the midpoint of the path from $x$ to $x'$ and hence that if $x<x'$ then $z$ 
is closer to $x'$ than it is to $x$ if and only if $z \in m\{x, x'\}$.  

\end{section}
\begin{section}{Theorem}
\begin{theorem}
The question: Is a given midpoints structure (or triples structure) realizable?  is NP complete.  
\end{theorem}
\begin{proof}
We will encode 3-SAT and then apply \ref{r>s} and \ref{s>r}.  
\end{proof}
\begin{corollary} 
Determining a compatible metric tree structure given a realizable  
midpoints structure (or triples structure) is NP hard.  
\end{corollary} 
\end{section}
\begin{section}{Encoding}
We will encode a case of 3-satisfiability in conjunctive normal form 
with $V$ variables and $C$ length 3 {\it or} clauses.  
The function $I$ names the variables appearing in a particular clause, 
while $\sigma$ indicates whether each variable appears with a {\it not}.  
We assume that all three variables appearing in any given clause are 
distinct.  

If $a\leq b \in {\bf Z}$ write $[a,b]$ for $\{a, a+1, \ldots ,b\}$.  
\begin{definition}
A {\bf case of 3-SAT} is a quadruple $(V,C,\nu,\sigma)$ with $V,C \in {\bf N}$, 
$\nu:[1,C] \times [0,2] \rightarrow [1,V]$ and 
$\sigma:[1,C] \times [0,2] \rightarrow \{-1,+1\}$ with $\nu(c,0)<\nu(c,1)<\nu(c,2)$ 
for every $c \in [1,C]$.  
The case $P=(V,C,\nu,\sigma)$ is said to be {\bf satisfiable by $h$} if 
$h: [1,V] \rightarrow \{-1,+1\}$ 
and for each $c \in [1,C]$ there is some $a \in [0,2]$ with $h(\nu(c,a))=\sigma(c,a)$.  
\end{definition}
\begin{example}\label{ex1} The case $(x \vee \overline  y \vee z)\wedge(w\vee x \vee y)$ 
is encoded with $V=4$, $C=2$, $\nu(1,1)=2$, $\sigma(1,1)=1$, ...
$$\{\nu(c,a)\}=\left(\begin{array}{ccc}
2&3&4\\1&2&3\end{array}\right) \hbox{  and  }
\{\sigma(c,a)\}=\left(\begin{array}{ccc}
+1&-1&+1\\+1&+1&+1\end{array}\right).$$
and is satisfiable by $11$ choices of signs one of which is 
$h(1)=-h(2)=-h(3)=h(4)=1$.  
\end{example}

The idea will be to construct a single midpoints structure for each 
case of 3-SAT.  This will be accomplished by using two constructions 
and combining them (one copy of the first and many of the second) 
using the following definition.  

\begin{definition}\label{cup}{Combining Midpoints Structures:}
\newline If $m:{X\choose 2} \rightarrow 2^X$ and $n:{Y\choose 2} \rightarrow 2^Y$ 
are midpoints structures on ordered sets $X$ and $Y$ respectively and $f:Y \rightarrow X$, 
then $m\cup_f n$ is a midpoints structure on $X \cup Y$ (with $X$ and $Y$ 
ordered as before and $x<y$ for every 
$x\in X$ and $y\in Y$) defined by $$\left. \begin{array}{rll}
(m\cup_f n)\{x,x'\}=&m\{x,x'\}\cup f^{-1}m\{x,x'\}& \hbox{ if } x,x' \in X\\
(m\cup_f n)\{x,y\}=&\{y\}& \hbox{ if } x \in X, y \in Y\\
(m\cup_f n)\{y,y'\}=&n\{y,y'\}& \hbox{ if } y,y' \in Y
\end{array} \right. $$
\end{definition}
If $m\cup_f n$ 
is realizable any realization restricts to realizations of $m$ and $n$ 
but if $m$ and $n$ are both realizable $m\cup_f n$ might or might not be.  
The idea for realizing $m\cup_f n$ 
when possible is to choose the distances for $m$ 
to be much smaller than those for $n$ and put the tree for $m$ into the 
middle of the tree for $n$ with leaves $Y$ attached to leaves 
$X$ according to the map $f$.  
\begin{example}
If $X=\{a,b,c\}$, $Y=\{A,B,C,D\}$, $f(A)=a$, $f(B)=b$, $f(C)=c$,  
$f(D)=c$  and the first two figures represent tree metrics realizing 
$m$ and $n$ respectively then the third figure represents a tree 
metric realizing $m \cup_f n$:
\newline\epsfig{file=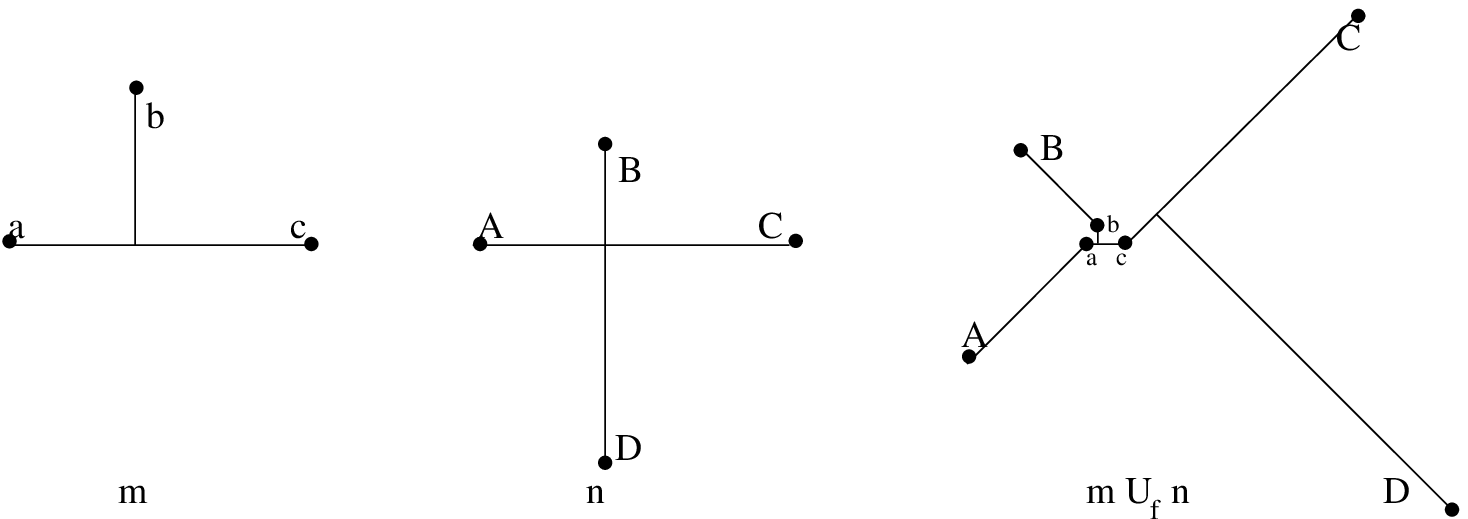}
\end{example}
For the remainder of this paper fix a case $P=(V,C,\nu,\sigma)$ of 3-SAT.  
We will construct in several steps a midpoints structure $m_P$ associated to $P$. 
The construction will be polynomial size and require polynomial time in the size of $P$ 
and $m_P$ will be realizable iff $P$ is satisfiable.  
\begin{subsection}{Variable Structure}\label{varms}
We start with a midpoints structure $m_0$ with two elements for 
each variable and four extra end elements.  
For each clause a copy of a second midpoints structure $m$ 
will be attached to $m_0$ using definition \ref{cup}.  
\newline 
The set of elements for $m_0$ is $X_0=\{x_{v,s}|v \in [0,V+1], s\in \{-1,1\}\}$ 
making a total of $2V+4$ elements.  
Order these elements with $x_{v',s'}<x_{v,s}$ if $v'<v$ or 
if $v'=v$ and $s'<s$.  
If $x_{v',s'}<x_{v,s} \in X_0$ then set 
$$m_0\{x_{v',s'},x_{v,s}\}:=
\left\{ \begin{array}{ll} 
\{x_{u,r} \in X_0|u\geq v, r \in \{1, -1\}\} & \hbox{ if } 2v'+s'=2v+s \hbox{ and }v'<v \\ 
\{x_{v,s}\} & \hbox{ if } 2v'+s'<2v+s \end{array} \right. .$$  
Note that $m_0$ is realizable with its midpoints tree geometry.  
For instance the vector $t\in {\bf R}^{2^{X_0}}$ with 
$t(\{x_{u,r}|u \geq v,r \in \{1,-1\}\})=1$
and $t(\{x_{v,s}\})=10^{v+{s+1 \over 2}}$ except that $t(\{x_{0,-1}\})=0$ 
and with $t(S)=0$ for all other $S \subseteq X_0$ 
is a realization of $m_0$.  
\end{subsection}

\begin{subsection}{Example:}\label{ex}
For example \ref{ex1} above $X_0=\{x_{0,-1}$, $x_{0,1}$, $x_{1,-1}$, $x_{1,1}$, $x_{2,-1}$, $x_{2,1}$, 
$x_{3,-1}$, $x_{3,1}$, $x_{4,-1}$, $x_{4,1}$, $x_{5,-1}$, $x_{5,1}\}$ 
as an ordered set, $m_0\{x_{0,1}, x_{2,1}\}=\{x_{2,1}\}$, 
$m_0\{x_{3,1}, x_{4,-1}\}=\{x_{4,1}$, $x_{4,-1}$, $x_{5,1}$, $x_{5,-1}\}$,  ...
and the figure represents the realization $t$ given above.  

\epsfig{file=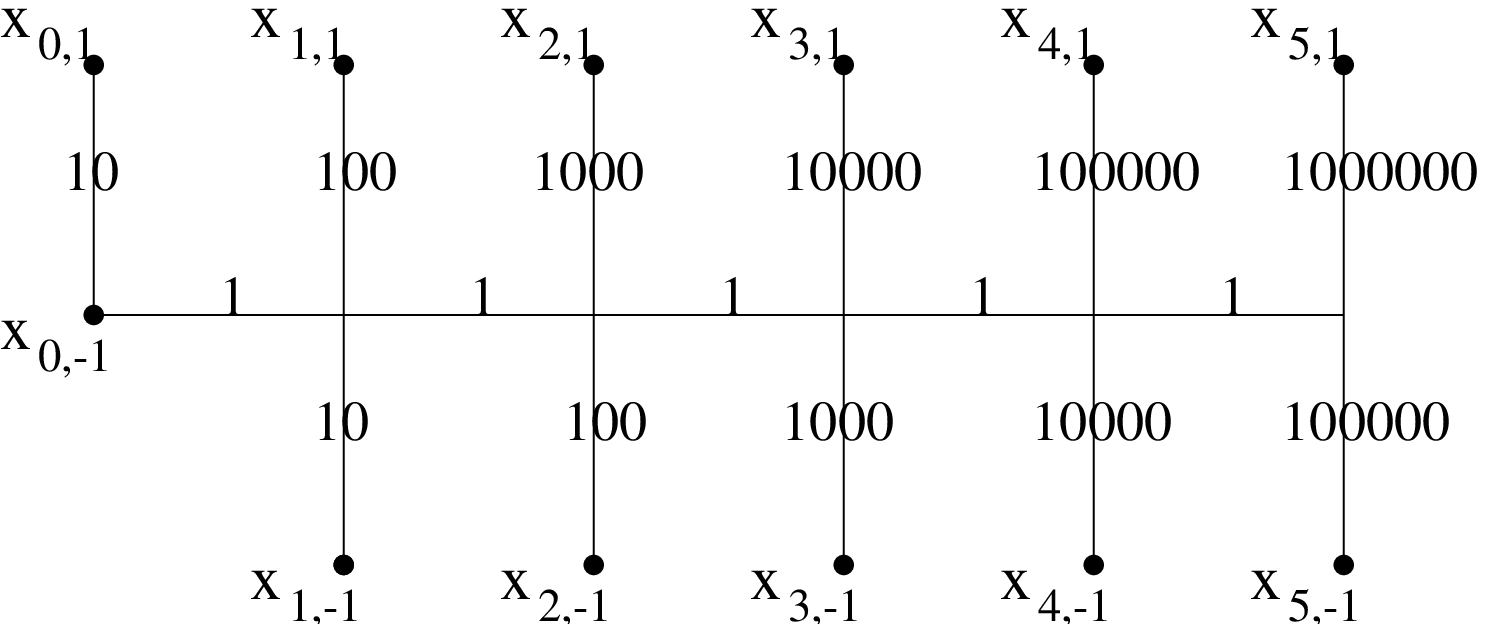, width=4in}
\end{subsection}
\vskip15pt\begin{subsection}{Encoding Signs:}\label{+-}
By lemma \ref{lem3} below every realization of $m_0$ must have every edge in 
the image of $m_0$ with a positive length and hence every realization must 
have a contraction to the midpoints tree.  Thus every realization 
involves pulling apart some (or none) of the $V+2$ degenerate vertices in this tree.  
For each degree four vertex there are $4$ ways to do this.  
\newline\epsfig{file=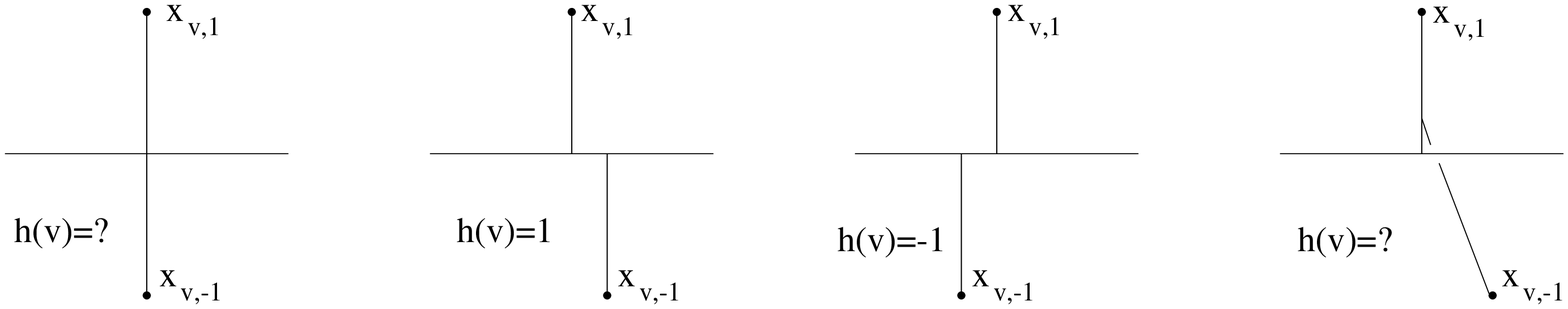, width=6in}
\vskip5pt
The idea of the construction will be to encode a choice of 
$h:[1,V]\rightarrow \{1,-1\}$ by the choice of this splitting 
with the second configuration corresponding to $h(v)=1$ 
and the third corresponding to $h(v)=-1$.  
The other two configurations will correspond to $h(v)$ being undetermined.  
We will add on a copy of the midpoints structure $m$ given below 
for each clause so that a splitting of the branches arises in a realization 
of the final midpoints structure if and only if the associated choice of 
$h:[1,V]\rightarrow \{1,-1\}$ satisfies the given case of 3-SAT.  
If a value $h(v)$ is undetermined it must satisfy the case for both values.  
\end{subsection}
\begin{subsection}{Clause Structures}\label{clams}
In this section we give the main construction.  
This is a midpoints structure $m$ with $48$ elements which has the property 
that after being added to $m_0$ there are at least three minimal sets of edges which can appear as 
a tree structure realizing the result.  This means that there is no realization
with the midpoints tree geometry since by lemma \ref{lem3} 
this geometry is a contraction of every geometry realizing a midpoints 
structure.  
This makes the sum of $m$ and $m_0$ a counterexample to Warnow's question.  
In this case there are three other edges 
so that every tree structure containing one of these three 
along with the edges 
containing midpoints is a realization of  $m$. 
We will use a copy of $m$ for each clause $c$ of the 
case $P$ of 3-SAT and denote it by $m_c$.  These will be combined with 
$m_0$ from above using definition 4.3 to get the final midpoints structure $m_P$.  
The clause is encoded by the choice of map with which to 
combine $m_c=m$ with $m_0$.  
The three special edges mentioned above will be positioned by the 
combining map so as to correspond to the three variable sign choices \ref{+-}
which will satisfy the $c$th three term {\it or} clause.  
The midpoints 
structure $m$ splits into $4$ isomorphic substructures 
each with $12$ elements (indexed by $p$ below) 
with $m$ applied to two elements with different $p$ indices being simply 
the singleton with larger $p$ so that in any realization the different $12$ 
element subtrees will have successively much longer leaf lengths.  
Each $12$ element midpoints substructure is simply a total order with 
a realization by its midpoints geometry which is simply a star with the 
ordering on leaf lengths given by the order below \ref{m}.  
The other $12$ element substructures differ only by a relabeling.  
The arrangement of the elements of these $12$ element substructures given 
by the combining maps to the 
underlying variable midpoints structure $m_0$ 
is shown in the figure after example \ref{example}.
The subtlety in $m$ arises from the fact 
that some of the midpoints between two elements of a $12$ element substructure 
are not singletons but rather doubletons 
involving an element from an adjacent substructure.  
These give $16$ more edges than in the star in the midpoints geometry 
of $m$ so that the midpoints geometry has one degree $32$ vertex adjacent 
to 16 leaves and $16$ degree $3$ vertices each adjacent to two more leaves.  
To get a counterexample to theorem 2 of \cite{g7} choose $P$ to be any 
non-satisfiable case of 3-SAT and consider $m_P$.  

The midpoints structure $m$ is the central construction of this paper 
and is given explicitly below.  

\begin{definition}\label{m}
$Y=[0,3]\times[0,3]\times[-1,1]$.
Define an involution $\mu$ on $Y$ by $\mu(p,q,e)=((p+e)_4,(q+e)_4,-e)$
where $0 \leq (r)_4\leq 3$ is the reduction of $r$ modulo $4$.  
Every midpoint will be either a singleton or a doubleton $\{y, \mu y\}$.  
Consider the total ordering on $[0,3]\times [-1,1]$ with 
$(0,0)<(1,0)<(0,1)<(1,-1)<(2,0)<(1,1)<(2,-1)<(3,0)<(3,1)<(0,-1)<(2,1)<(3,-1)$
Totally order $Y$ by setting $(p',q',e')<(p,q,e)$ if $p'<p$ or if $p'=p$ and 
$(q',e')<(q,e)$ above. 
Define $m$ so that for every $(p',q',e')<(p,q,e)\in Y$ we have 
$m\{(p',q',e'),(p,q,e)\}=\{(p,q,e)\}$ if 
$p'<p$ or if [$p'=p$ and $(q',e')<(q,0)$ above] or if 
[$p'=p$, $q'=q$, $e'=0$ and $e=1$] and otherwise 
$m\{(p',q',e'),(p,q,e)\}=\{(p,q,e),\mu(p,q,e)\}$.  
\end{definition}

Recall that $P=(V,C,\nu,\sigma)$ is an instance of 3-SAT.  
For each $c \in [1,C]$ define $m_c$ to be the midpoints structure on 
$X_c=\{x_{c,p,q,e}|p,q \in [0,3], e\in [-1,1]\}$ isomorphic to $m$ above. 
Extend the total order and the involution $\mu$ from $Y$ to $X_c$.  
(Explicitly take $\phi :X_c \rightarrow Y$ to be $\phi x_{c,p,q,e}=(p,q,e)$ 
and set $m_c\{x,x'\}=\phi^{-1}m\{\phi x,\phi x'\}$, 
$x<x'$ if $\phi x<\phi x'$ and $\mu x=\phi^{-1} \mu \phi x$ 
for every $x, x' \in X_c$.)  
The following maps will be used to combine the $m_c$ with $m_0$ using 
\ref{cup}.  
Define $X_{\leq 0}:=X_0$ and $X_{\leq c}:=X_{\leq (c-1)}\cup X_c$ for every 
$c \in [1,C]$ and 
$X:=X_P:=X_{\leq C}$.  
Define $f:X \rightarrow X_0$ by setting $f|_{X_0}$ to be the identity and
$$\left.\begin{array}{rll} 
fx_{c,p,0,0}:=& x_{\nu(c,0),-\sigma(c,0)} & \\ 
fx_{c,p,1,0}:=& x_{\nu(c,1),-\sigma(c,1)} & \\ 
fx_{c,p,2,0}:=& x_{\nu(c,2),-\sigma(c,2)} & \\ 
fx_{c,p,3,0}:=& x_{0, 1} & \\
fx_{c,p,0,-1}:=fx_{c,p,3,1}:=& x_{\nu(c,0),\sigma(c,0)} & \\
fx_{c,p,1,-1}:=fx_{c,p,0,1}:=& x_{\nu(c,1),\sigma(c,1)} & \\
fx_{c,p,2,-1}:=fx_{c,p,1,1}:=& x_{\nu(c,2),\sigma(c,2)} & \\
fx_{c,p,3,-1}:=fx_{c,p,2,1}:=& x_{V+1, 1} & \end{array} \right. $$
\vskip0pt
Note that $f$ is invariant under the involution $\mu:X \rightarrow X$.  
We combine sequentially to get the final construction of the midpoints structure $m_P$.  
\begin{definition}
If $P=(V,C,\nu,\sigma)$ is a case of 3-SAT then take $m_{\leq 0}:=m_0$ from above and 
for every $c \in [1,C]$ take $m_c$ and $f$ from above and
define $m_{\leq c}:=m_{\leq (c-1)}\cup_{f|_{X_c}}m_c$ a 
midpoints structure on $X_{\leq c}$.  
Finally take $m_P:=m_{\leq C}$ the midpoints structure on $X$.  
\end{definition}
\end{subsection}
%\newpage
\begin{example}\label{example}
For the example \ref{ex1} above we have $C=2$, 
$X_1=\{x_{1,0,0,-1}$, $x_{1,0,0,0}$, \ldots, $x_{1,3,3,1}\}$ and 
$X=X_0\cup X_1\cup X_2$ with $108$ elements.  
The map $f:X \rightarrow X_0$ has for instance $f^{-1}(x_{0,1})
=\{x_{0,1}, x_{1,0,3,0}$, $x_{1,1,3,0}$, $x_{1,2,3,0}$, $x_{1,3,3,0}$, 
$x_{2,0,3,0}$, $x_{2,1,3,0}$, $x_{2,2,3,0}$, $x_{2,3,3,0}\}$, 
$f^{-1}(x_{1,1})=\{x_{1,1}$, $x_{2,0,0,-1}$, $x_{2,0,3,1}$, $x_{2,1,0,-1}$, 
$x_{2,1,3,1}$, $x_{2,2,0,-1}$, $x_{2,2,3,1}$, $x_{2,3,0,-1}$, $x_{2,3,3,1}\}$, 
and $f^{-1}(x_{1,-1})=\{x_{1,-1}$, $x_{2,0,0,0}$, $x_{2,1,0,0}$, 
$x_{2,2,0,0}$, $x_{2,3,0,0}\}$.
Below is part of a realization of $m_P$ based on the 
solution $h$ from example \ref{ex1}.  
Only $26$ of the $108$ elements are shown 
(the $12$ in $X_0$, $13$ of those from $X_1$ and one from $X_2$).  The lengths 
given are those of $t_h$ given in \ref{s>r}.
The adjustment $t'$ to $t$ given in \ref{s>r} is the difference between the 
labeled leaf lengths and the nearest multiple of $10^{12}$ 
(the position in the numbers below marked with a semicolon); $t'$ does not adjust the interior edges.  
\vskip10pt\hskip-15pt
\epsfig{file=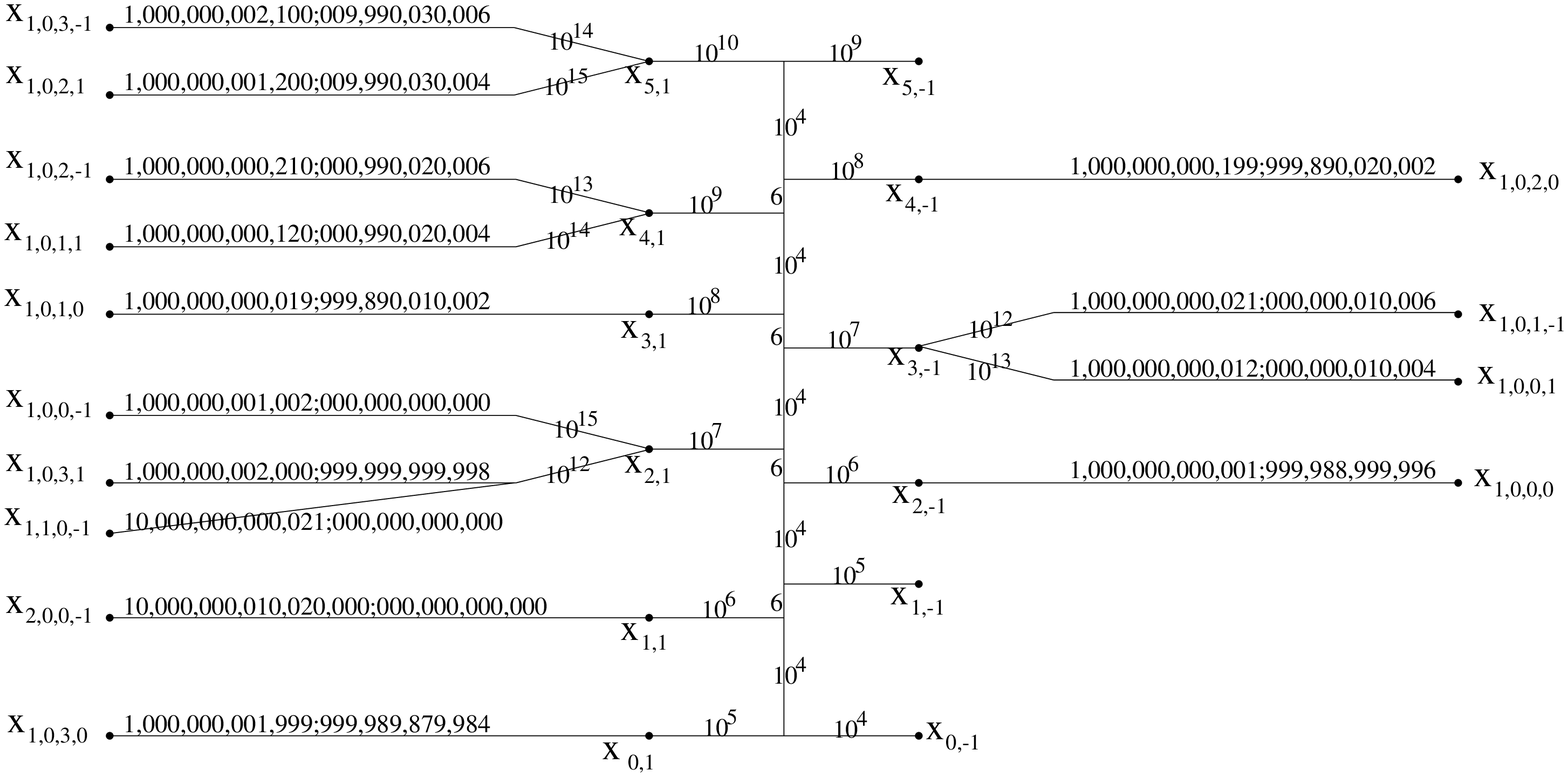,width=6.5in}
\end{example}
\end{section}
%\newpage
\begin{section}{Realizable implies satisfiable}
\begin{proposition}\label{r>s}
If $P$ is a case of 3-SAT and $m_P$ is realizable 
then $P$ is satisfiable.  
\end{proposition}

Recall definition \ref{int}.  
\begin{lemma}\label{lem1}If $t$ is a realization of  some midpoints structure on $X$ and $S,T,U \subseteq X$  with 
$t(T)t(S)t(U)\not=0$ and $T \in [S,U]$ (that is $S:T:U$) 
then $t_{(S,U)}=t_{(S,T)}+2t(T)+t_{(T,U)}$. 
\end{lemma}
\begin{proof}
This follows geometrically from the fact that $t_{(S,U)}$ is twice the distance 
in the tree between the ends of the edges $S$ and $U$ which are closest 
to each other, and $T \in [S,U]$ with $t(T)>0$ iff $T$ is an edge 
on the unique geodesic in the tree $t$ between $S$ and $U$.  
\end{proof}
\epsfig{file=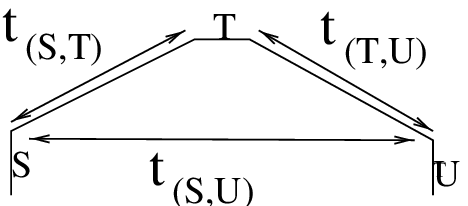}

\begin{lemma}\label{lem2} If $t$ is a realization of  some midpoints structure on $X$ and 
$U_0 \cup U_1 \cup U_2 \cup U_3 =X$ is a partition of $X$ 
with $t(U_i)>0$ for every $i$ then $t(U_j\cup U_3)>0$ 
for at most one $j \in [0,2]$.  Further, $t_{( U_j, U_k)}=2t(U_j \cup U_k)$
so the sum 
$\tau=\Sigma_{j\in [0,3]} (-1)^jt_{( U_j, U_{(j+1)_4})}=4t(U_0 \cup U_3)-4t(U_2 \cup U_3)$ 
and hence is positive 
iff $t(U_0\cup U_3)>0$ (so $j=0$) and negative iff 
$t(U_2\cup U_3)>0$ (so $j=2$).  If $j=1$ this sum is $0$.  
\end{lemma}
\begin{proof}
For all $e,f \in \{1,c\}$ we have $(U_j\cup U_3)^e\cap(U_k\cup U_3)^f\in \{U_i\}_{i \in [0,3]}$ 
and since each $U_i$ is nonempty we get that $t(U_j\cup U_3)t(U_k\cup U_3)=0$ 
by definition \ref{real}.  For the second statement if $T \in (U_j, U_k)$ with $t(T)>0$
then using 
definition \ref{int} and replacing $T$ with $T^c$ if necessary we get 
$U_j \subseteq T \subseteq U_k^c$.  Since we are taking the open interval 
$T \not=U_j$ and $T\not=U_k^c$.  Choose $U_i\not\in \{U_j, U_k\}$ with 
$T \cap U_i \not= \emptyset$.  By definition \ref{real} there is some choice of 
signs $e,f \in \{1,c\}$ with $T^e \cap U_i^f=\emptyset$ so $U_i \subseteq T$ and 
hence $T=U_i \cup U_j$.  
\end{proof}
\epsfig{file=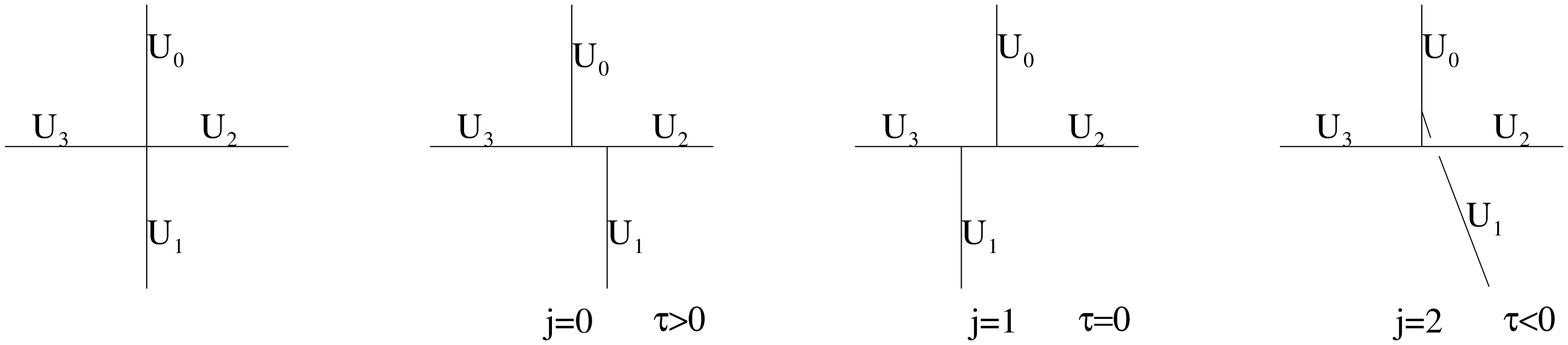, width=6in}
\begin{lemma}\label{lem3} For any realization $t$ of a midpoints structure $m$ on $X$ 
and any $x \not= y \in X$ we have $t(m\{x,y\})>0$, so 
the edge containing any midpoint must have positive length and hence 
every tree geometry for which there is a realization of $m$ contains the 
midpoints geometry as a contraction.  
\end{lemma}
\begin{proof}
By (3) in the definition of realization \ref{real} we have 
$t_{[x,m\{x,y\}]}>t_{[y,m\{x,y\})}$ and similarly, switching 
the roles of $x$ and $y$ we get $t_{[y,m\{x,y\}]}>t_{[x,m\{x,y\})}$.  
Adding these and canceling gives $4t(m\{x,y\})>0$.  
\end{proof}
\begin{proof}{\it of Proposition:}
First we will define a choice of signs $h_t:[V] \rightarrow \{1,-1\}$ 
given a realization $t$ of $m_P$.  
We will then check that $h_t$ satisfies $P$.  

Assume that $t \in {\bf R}^{2^X}$ is a realization of $m_P$.  
For convenience we give new names to some of the elements $A \subseteq X$ in the 
image of $m_P$ and their complements.  Write 
$$\left. \begin{array}{rcl}
&\{x_{c,p,q,e}\}&=m_P\{x_{c,p,q,e}, x_{v,s}\}, \\
\{x_{c,p,q,-1}\}^\mu :=\{\mu x_{c,p,q,-1}\}^\mu :=&\{x_{c,p,q,-1}, \mu x_{c,p,q,-1}\}&=
 m_P\{x_{c,p,q,-1}, x_{c,p,(q-1)_4,+1}\}, \\
A_{v,s}:=&f^{-1}\{x_{v,s}\}&=m_P\{x_{v,s}, x_{v-1, s}\}, \\
A_{>v}:=(A_{<v-1})^c:=&f^{-1}\{x_{v',s}|v'>v, s\in \{1,-1\}\}&=m_P\{x_{v,1}, x_{v+1,-1}\}. \\
\end{array} \right.$$
For every $v$ and $s$ the set 
$\{A_{v,s}, A_{v,-s}, A_{>v}, A_{<v}\}$ partitions $X$.  

Now by lemma \ref{lem3} we get that for every midpoint above 
$t(m_P\{x,y\})>0$.  
Thus we can apply lemma \ref{lem2} to each of the (ordered) partitions above.  
Write $\tau_{v,s}=t_{(A_{v,s}, A_{v,-s})}-
t_{(A_{v,-s}, A_{>v})}+t_{(A_{>v}, A_{<v})}-t_{(A_{<v}, A_{v,s})}$ 
which by \ref{lem2} 
will have $\tau_{v,s}>0$ (and $\tau_{v,-s}=0$) if and only if  
$t(A_{v,s} \cup A_{<v})>0$.  
Since replacing $s$ with $-s$ simply switches the 
order of the first two elements of the partition, 
\ref{lem2} also gives that $\tau_{v,-s}>0$ 
(and $\tau_{v,s}=0$) iff $t(A_{v,-s} \cup A_{<v})>0$.  
In particular these situations are mutually exclusive.  
\begin{definition}
Define $h_t(v)=s$ if $\tau_{v,s}>0$ and choose $h_t(v)$ arbitrarily if 
$\tau_{v,1}=\tau_{v,-1}\leq 0$.  
\end{definition}

We will show that $h_t$ satisfies $P$.  
Fix a clause $c \in [1,C]$.  
For each fixed $p\in[0,3]$ order lexicographically the 12 elements 
$\{x_{c,p,q,e}|q \in [0,3], e \in [-1,1]\}\subseteq X_c$ 
(so $x_{c,p,0,-1}<$ $x_{c,p,0,0}<$ $x_{c,p,0,1}<$ $x_{c,p,1,-1}<$ $x_{c,p,1,0}<$
$x_{c,p,1,1}<$ $x_{c,p,2,-1}<$ $x_{c,p,2,0}<$ $x_{c,p,2,1}<$ $x_{c,p,3,-1}<$
$x_{c,p,3,0}<$ $x_{c,p,3,1}$)
and consider the 12 inequalities $0<t_{[\{x\},m\{x,y\}]}-t_{[\{y\},m\{x,y\})}$ 
(from part (3) of \ref{real}) 
obtained by taking $y<x$ to be adjacent in the above order (11 cases) 
or else $x=x_{c,p,3,1}$ the last element and $y=x_{c,p,0,-1}$ the first.  
Add the 48 inequalities obtained by taking these 12 inequalities for all 
choices of $p\in [0,3]$ to obtain:  

$$ 0< \sum_{p,q \in [0,3]}\big{(}
 t_{[\{x_{c,p,q,0}\}, \{x_{c,p,q,-1}\}^\mu]}
-t_{[\{x_{c,p,q,-1}\}, \{x_{c,p,q,-1}\}^\mu)}
$$
$$
+t_{[\{x_{c,p,q,1}\}, \{x_{c,p,q,1}\}]}
-t_{[\{x_{c,p,q,0}\}, \{x_{c,p,q,1}\})}
+t_{[\{x_{c,p,q,-1}\}, \{x_{c,p,q,-1}\}^\mu]}
-t_{[\{x_{c,p,(q-1)_4,1}\}, \{x_{c,p,q,-1}\}^\mu)}\big{)}
$$

For every $c\in[1,C]$, $p\in[0,3]$ and $q\in[0,3]$ 
for the first line, $a\in[0,2]$ for the second and $b\in[0,1]$ for the fourth
we have:
$$\left. \begin{array}{c}
\{x_{c,p,q,1}\}:\{x_{c,p,q,1}\}^\mu:\{x_{c,p,(q+1)_4,-1}\}^\mu:\{x_{c,p,(q+1)_4,-1}\} \\
\{x_{c,p,a,0}\}:A_{\nu(c,a),-\sigma(c,a)}:A_{\nu(c,a),\sigma(c,a)}:\{x_{c,p,a,-1}\}^\mu:\{x_{c,p,a,-1}\} \\
\{x_{c,p,3,0}\}:A_{0,1}:A_{<\nu(c,0)}:A_{>\nu(c,0)}:A_{<\nu(c,1)}:A_{>\nu(c,1)}:A_{<\nu(c,2)}:A_{>\nu(c,2)}:A_{V+1,1}:\{x_{c,p,3,-1}\}^\mu:\{x_{c,p,3,-1}\} \\
\{x_{c,p,b,0}\}:A_{\nu(c,b),-\sigma(c,b)}:A_{>\nu(c,b)}:A_{<\nu(c,b+1)}:A_{\nu(c,b+1),\sigma(c,b+1)}:\{x_{c,p,b,1}\}^\mu:\{x_{c,p,b,1}\} \\
\{x_{c,p,2,0}\}:A_{\nu(c,2),-\sigma(c,2)}:A_{>\nu(c,2)}:A_{V+1,1}:\{x_{c,p,2,1}\}^\mu:\{x_{c,p,2,1}\} \\
\{x_{c,p,3,0}\}:A_{0,1}:A_{<\nu(c,0)}:A_{\nu(c,0),\sigma(c,0)}:\{x_{c,p,3,1}\}^\mu:\{x_{c,p,3,1}\} \\
\end{array} \right. $$
Using lemma \ref{lem1} and the above betweenness relations 
to substitute into the previous inequality gives (after some cancellation): 

$$
0< -\sum_{p,q\in[0,3]}(2t_{(\{x_{c,p,q,1}\}, \{x_{c,p,q,1}\}^\mu)}+
t_{(\{x_{c,p,(q-1)_4,-1}\}^\mu, \{x_{c,p,q,-1}\}^\mu)})
$$
$$
+4\sum_{a\in[0,2]}(t_{(A_{\nu(c,a),\sigma(c,a)}, A_{\nu(c,a),-\sigma(c,a)})}
                 -t_{(A_{\nu(c,a),-\sigma(c,a)}, A_{>\nu(c,a)})}
                 +t_{(A_{>\nu(c,a)}, A_{<\nu(c,a)})}
                 -t_{(A_{<\nu(c,a)}, A_{\nu(c,a),\sigma(c,a)})})
$$
$$
\leq 4\sum_{a\in[0,2]}\tau_{\nu(c,a),\sigma(c,a)}
$$ 
Thus there is some $a \in [0,2]$ for which 
$\tau_{\nu(c,a),\sigma(c,a)}>0$ and hence $h_t(\nu(c,a))=\sigma(c,a)$.  
\end{proof}
\end{section}
\begin{section}{Satisfiable implies realizable}
\begin{proposition}\label{s>r}
If $P$ is a satisfiable case of 3-SAT then the 
midpoints structure $m_P$ on $X_P$ is realizable.  
\end{proposition}
\begin{proof}
Assume that $h:[1,V] \rightarrow \{-1,1\}$ satisfies $P$.  
Construct $t_h \in {\bf R}^{2^X}$ realizing $m_P$ by starting with 
$t \in {\bf R}^{2^X}$ which almost realizes $m_P$ and then perturbing 
the values on the leaf edges by $t'$ to get $t_h$.  
Start with  
$$\left. \begin{array}{rl}
t(A_{>v}\cup A_{v,h(v)})= & 6 \hbox{ \hskip10pt for all }v \in [1,V], \\
t(A_{>v})= & 10^{V} \hbox{  \hskip10pt for all }v \in [0,V], \\
t(A_{v,s})= & 10^{V+v+{s+1 \over 2}} \hbox{ \hskip10pt  for all }v \in [0,V+1] 
\hbox{ and }s \in \{-1,1\}, \\
t(\{x_{c,p,q,e}\}^\mu)= & 10^{2V+4c+(p+q+e)_4} 
\hbox{  \hskip10pt for all }c \in [1,C]\hbox{, }p,q \in [0,3]\hbox{ and }e \in \{-1,1\},\\
t(\{x_{c,p,q,e}\})= & 10^{2V+4C+4+4c+p} 
+2\times 10^{2V+4c+(p+q)_4}+e^2 \times 10^{2V+4c+(p+q+e)_4} \\
&\hbox{ for all }c \in [1,C]\hbox{, }p,q \in [0,3]\hbox{ and } e \in [-1,1] \\
 \end{array} \right.$$
and $t(S)=0$ for all other $S \subseteq X$. 
Note that the only edges of $t$ (and $t_h$) not containing midpoints are 
$A_{0,-1}$ and those in the first line above with length $6$.  
The vector $t$ fails to realize $m_P$ only for the midpoints 
of lexicographically adjacent pairs of elements $x_{c,p,q,e}$ of $X$ with the 
same $c$ and $p$ coordinates.  This is corrected by the slight 
perturbation below.  
Denote the number of agreements by $n_h(c)=|\{a \in [0,2]| 
h(\nu(c,a))=\sigma(c,a)\}|$.  

Correction:  $$\left. \begin{array}{rl}
t'(\{x_{c,p,0,-1}\})= &0, \\
t'(\{x_{c,p,q,0}\})= &t'(\{x_{c,p,q,-1}\})-u_{[fx_{c,p,q,-1},fx_{c,p,q,0}]}+n_h(c), \\
t'(\{x_{c,p,q,1}\})= &t'(\{x_{c,p,q,0}\})+u_{[fx_{c,p,q,0},fx_{c,p,q,1}]}+n_h(c), \\
t'(\{x_{c,p,q,-1}\})= &t'(\{x_{c,p,q-1,1}\})+n_h(c) \hbox{ if } 
q \not= 0\\
 \end{array} \right.$$ 
and $t'(S)=0$ for all other $S\subseteq X$.  
Finally set $t_h(S)= t(S)+t'(S)$ for every $S \subseteq X$.  

It is now straightforward to check that $t_h$ is a realization of $m_P$.  

%In particular:  
%The correction terms have $|t'(\{x_{c,p,q,e}\})|<
%47\times (\hbox{max}\{u_{[x_{i,e}, x_{i',e'}]}\}+3)<
%47\times 2(10^{2V+2}+10^{2V+1}+(V+1)(10^V+6))<10^{2V+4}$.  
%Thus in checking the inequality from \ref{real} part (3) for $t$ 
%it will suffice to check that if $x \not= y$ then 
%$u_{[\{x\}, m\{x,y\}]} > u_{[\{y\}, m\{x,y\})} +10^{2V+4}$.  
%This will suffice for many of the pairs $x,y$.  
%If $x=x_{c,p,q,s}$ and $y=x_{i,e}$ then 
%$m\{x, y\}=\{x\}$ so $u_{[\{x\}, m\{x,y\}]} =
%2\times(10^{2V+4C+4+4c+p}+2\times 10^{2V+4c+(p+q)_4}+10^{2V+4c+(p+q+s)_4})$ 
%while $ u_{[\{y\}, m\{x,y\})}=u_{[\{y\}, \{fx\}]}<10^{2V+4}$ as above.  
%Similarly if $y=x_{c,p,q,s}$ and $x=x_{i,e}$ the inequality becomes 
%$\cdots>0$ which will hold trivially.  

\end{proof}
\end{section}

\bibliographystyle{amsplain}
\bibliography{biblio}
\end{document}